\documentclass[12pt]{amsart}

\usepackage{amsmath}
\usepackage{amssymb}

\newtheorem{theo}{Theorem}[section]
\newtheorem{definition}[theo]{Definition}
\newtheorem{lemma}[theo]{Lemma}
\newtheorem{prop}[theo]{Proposition}
\newtheorem{cor}[theo]{Corollary}
\theoremstyle{remark}
\newtheorem{example}[theo]{Example}
\newtheorem{remark}[theo]{Remark}

\def\N{{\mathbb N}}

\def\Q{{\mathbb Q}}
\def\T{{\mathbb T}}
\def\cT{{\mathcal T}}
\def\Nk{{{\mathbb N}^k}}

\def\Z{{\mathbb Z}}

\newcommand{\clsp}{\overline{\operatorname{span}}}

\newcommand{\End}{\operatorname{End}}
\numberwithin{equation}{section}

\begin{document}

\title[Crossed products via transfer operators]{Crossed 
products by abelian semigroups via transfer operators}
\author[N. S. Larsen]{Nadia S. Larsen}
\address{Department of Mathematics, University of Oslo,
P.O. Box 1053 Blindern, N-0316 Oslo, Norway.}
\email{nadiasl@math.uio.no}
\thanks{Research partly supported by The Carlsberg Foundation and the 
Centre for Advanced Study, Oslo.}
\date{February 15, 2005}
\begin{abstract}
We propose a generalisation of Exel's crossed product by a single endomorphism 
and a transfer operator to the case of actions of abelian semigroups of 
endomorphisms and associated transfer operators. The motivating example for 
our definition yields new crossed products, not obviously covered by 
familiar theory. Our technical machinery 
builds on Fowler's theory  of Toeplitz and Cuntz-Pimsner algebras of 
discrete product systems of Hilbert bimodules, which we need to expand to cover  
a natural notion of relative Cuntz-Pimsner algebras of product systems. 
\end{abstract}
\maketitle

The crossed product of a $C^*$-algebra $A$ by 
an endomorphism $\alpha$ was first constructed by Cuntz \cite{C} as a corner in an 
ordinary group crossed product. Since then, other authors have 
proposed definitions of such a crossed 
product mainly in terms of universal properties of representations 
of $A$ and the unital, additive semigroup $\N$,  
cf. \cite{P, S, ALNR, M}. The construction introduced by Exel  
\cite{E} adds to the pair $(A, \alpha)$ 
a new ingredient, namely a \emph{transfer operator}, 
which is a positive continuous linear map $L:A\to A$  satisfying 
$L(\alpha(a)b)=aL(b)$ for all $a, b\in A$.  
The resulting crossed 
product $A\rtimes_{\alpha, L} \N$ generalises the 
previous constructions 
from \cite{P, S, ALNR, M}, see also the more recent \cite{ER}.

In Stacey's theory of crossed products by single  endomorphisms \cite{S}, the
extension to other semigroups brought in a large number of new and
interesting examples and applications, see for example \cite{A, ALNR, LR1, LR2}.

Here we seek an extension of Exel's theory, and we are motivated by the following 
situation. For a compact abelian group $\Gamma$ there is always an  
action $\beta:\N^* \to End(C(\Gamma))$ of the multiplicative semigroup $\N^*$ 
given by $\beta_n(f)(s)=f(s^n)$. What we notice is that under certain finiteness 
conditions there is an action of transfer operators $L_n$ for $\beta_n$, see  
Proposition~\ref{beta_L_when_injective}, and 
then we wish to form and analyse a crossed product $C(\Gamma)\rtimes_{\beta, 
L}\N^*$. It is our goal here to adapt Exel's 
definition to the case of actions of abelian semigroups of endomorphisms and 
corresponding transfer operators. As applications, 
we take a closer look at the motivating example. Our construction allows us 
to view 
in a unified way the crossed product by two rather different types of actions, 
see Theorem~\ref{swapping_Exel_cp}, and for the motivating example we use this 
to show 
that in some cases the crossed product is covered by existing theory. In  
other cases modelled by our example the resulting crossed product is  
new, and therefore of interest for further analysis.

 Our approach is based 
on Fowler's theory of discrete product systems of Hilbert bimodules 
over a semigroup \cite{F}, and is closely related to Brownlowe and Raeburn's 
work \cite{BR}, where they show that Exel's  
crossed product by a single endomorphism is a relative Cuntz-Pimsner 
algebra as in \cite{FMR, MS}. What we find is that developing a theory 
of crossed products by a semigroup of endomorphisms with 
associated transfer operators goes hand in hand with the need to 
introduce and study relative Cuntz-Pimsner algebras of product systems of 
Hilbert bimodules. Here we initiate this study along the lines 
of \cite{F} (see also \cite{Br2}).

The paper is organised as follows: in section $1$ we consider a semigroup 
dynamical system $(A, S, \alpha)$ consisting of a $C^*$-algebra $A$ and an action 
of a unital semigroup $S$ by endomorphisms of $A$, together with an action $L$ of 
transfer operators for $\alpha$; we associate to $(A, S, \alpha, L)$ a Toeplitz 
crossed product, which is a universal  $C^*$-algebra for certain covariant 
representations. As in \cite{E}, we define the crossed 
product $A\rtimes_{\alpha, L}S$ to be a quotient of this Toeplitz crossed 
product by an ideal of redundancies. In section $2$ we construct a product system $X$ 
of Hilbert bimodules over $S$ whose Toeplitz algebra \cite{F} is a Toeplitz 
crossed product for $(A, S, \alpha, L)$.

In section $3$ we identify $A
\rtimes_{\alpha, L}S$ as a natural analogue of a 
relative Cuntz-Pimsner algebra of the Toeplitz algebra of the product system $X$. 
In the case when $L$ is implemented by an action of injective 
endomorphisms with hereditary range we show in 
Proposition~\ref{inj_endom_hereditary_range} 
that the new crossed product 
$A\rtimes_{\alpha, L}S$ is the same as $A\rtimes_\alpha S$ from \cite{M, LR1, 
Lar}. 
Specialising to the semigroup $S=\N$ this result says that the relative Cuntz-Pimsner 
algebra $A\rtimes_{\alpha, L}\N$ of the product system $X$ over $\N$ is $A
\rtimes_\alpha \N$, 
which by \cite[Theorem 4.7]{E} is isomorphic to Exel's $A\rtimes_{\alpha, 
L}\N$. Since  
\cite[Proposition 3.10]{BR} shows that Exel's $A\rtimes_{\alpha, L}\N$
is isomorphic to the relative Cuntz-Pimsner algebra of a bimodule and an ideal 
arising from $A$, we are guided to predict that 
 our $A\rtimes_{\alpha, L}\N$, seen 
as a relative Cuntz-Pimsner algebra of a product system $X$ over $\N$, is isomorphic 
to the relative Cuntz-Pimsner algebra corresponding to the bimodule over the fibre 
at $n=1$ and its corresponding ideal. Indeed, we prove this as an application of the 
identification $A\rtimes_{\alpha, L}\N \cong A\rtimes_\alpha \N$, see 
Proposition~\ref{rel_CP_as_first_CP}. This collapsing of 
a quotient of the Toeplitz algebra of a 
product system $X$ over the semigroup $\N$ to the algebra over the fibre at $n=1$   
was noticed by Fowler in the case of the Cuntz-Pimsner algebra of a product 
system, see \cite[Proposition 2.11]{F}. Fowler's Proposition 2.11 has two sets of 
hypotheses, asking that the left action in the fibres is isometric or by 
compacts. We are grateful to Nathan Brownlowe for providing us with a copy of 
\cite{Br2}, where an extension of \cite[Proposition 2.11]{F} is obtained for relative 
Cuntz-Pimsner algebras with the first set of hypotheses and some further conditions. 
Our situation of the specific product system $X$ 
constructed from an action of $\N$ by injective 
endomorphisms with hereditary range can be regarded as an extension to relative 
Cuntz-Pimsner algebras of 
\cite[Proposition 2.11]{F} under an appropriate variant of the second set of hypotheses.

Section $4$ is devoted to our motivating class of examples. The general 
form  of these semigroup actions of $\N^*$ by endomorphisms with associated transfer 
operators of $C(\Gamma)$, where $\Gamma$ is a compact abelian group, is described in 
Proposition~\ref{beta_L_when_injective}. The class of examples of crossed products 
by endomorphisms and transfer operators arising from the actions in 
Proposition~\ref{all_phin_surj} is a new construction, not obviously covered 
by familiar theory. However, another large class of examples here, which is described 
in Proposition~\ref{general_form_varphi_inj}, recovers semigroup crossed products in the 
sense of \cite{M, ALNR,LR1} 
which model Hecke $C^*$-algebras, cf. \cite{LR2, LarR}.

In extending Exel's $A\rtimes_{\alpha,L}\N$ 
corresponding to 
$(A, \N, \alpha, L)$, we chose to look at abelian semigroups first of all 
because the examples we have here involve the semigroups $\Nk$ for 
$k\in \{1, 2,\dots \}\cup \{\infty\}$, and secondly, because in order 
to construct 
product systems of Hilbert bimodules over non-abelian semigroups as in  
\cite{F} we would have to impose technical conditions on our actions. 
We expect that some of our constructions could be extended to cover 
the case of Ore semigroups. 

We thank Iain Raeburn for many useful discussions and suggestions, 
Nathan Brownlowe for suggestions and for providing a copy of \cite{Br2}, and 
Magnus Landstad and Ola Bratteli for the invitation to attend the  
programme on ``Non-commutative Phenomena in Mathematics and Theoretical Physics'' 
at the Centre for Advanced Study, Oslo, 
during winter 2001, where the initial motivation for this work came about.

\section{The dynamical system $(A, S, \alpha, L)$}\label{the_system}

Let $A$ be a (not necessarily unital) $C^*$-algebra, and $S$ an 
abelian semigroup with identity $0$. Suppose that $\alpha:S \to \End(A)$  
is an action, i.e. a semigroup homomorphism, and assume that each $\alpha_s$ 
is \emph{extendible}, in the sense that it extends uniquely to an endomorphism 
$\overline{\alpha_s}$ of the multiplier algebra  
$M(A)$ such that $\overline{\alpha_s}(1_{M(A)})=\lim 
\alpha_s(u_\lambda)$ for  some (and hence every) approximate unit 
$(u_\lambda)$ in $A$ and all $s\in S$ (see \cite{A}).

We also assume that there is an action $L$ of $S$ by continuous, linear, 
positive maps 
$L_s:A\to A$ which admit linear continuous extensions $\overline{L_s}:M(A)\to 
M(A)$ that satisfy the \emph{transfer operator} identity
\begin{equation}
L_s(\alpha_s(a)u)=a\overline{L_s}(u)
\label{def_of_transfer_op} 
\end{equation}
for all $a\in A, u\in M(A), s\in S$. 
Note then that $L_s(u\alpha_s(a))=\overline{L_s}(u)a$ as well. 

A representation of $S$ in a 
$C^*$-algebra $B$ is a map $V:S\to M(B)$ such that 
$V_{st}=V_sV_t$ for all $s, t\in S$ and $V_0=1_{M(B)}$. If 
each $V_s$ is an isometry, we say $V$ is an isometric representation of $S$ 
in $B$. A non-degenerate representation 
$\pi:A\to B(H)$ and a representation $V:S \to B(H)$ form a 
\emph{Toeplitz covariant representation}, see \cite[Definition 3.1]{BR}, if they 
satisfy  the following conditions (inspired 
in relations (i)-(ii) of \cite[Definition 3.1]{E}):
\begin{enumerate}
\item[(TC1)]$V_s\pi(a)=\pi(\alpha_s(a))V_s$, and
\item[(TC2)]$V_s^*\pi(a)V_s=\pi(L_s(a))$, for all $a\in A$ and $s\in S$.
\end{enumerate}

\begin{definition} A Toeplitz crossed product for 
$(A, S, \alpha, L)$ is a triple $(B, i_A, i_S)$, where $i_A:A\to B$ is 
a non-degenerate homomorphism and $i_S:S\to M(B)$ is a representation 
satisfying the conditions
\begin{enumerate} 
\item[(i)]
if $\sigma$ is a non-degenerate representation of $B$, then $(\sigma\circ 
i_A, \overline{\sigma}\circ i_S)$ is a Toeplitz covariant representation;
\item[(ii)] 
for every Toeplitz covariant representation $(\pi, V)$, there is a  
representation $\pi\times V$ of $B$ such that $(\pi \times V)\circ i_A 
=\pi$ and $\overline{\pi \times V}\circ i_S=V$;
\item[(iii)]
$B$ is generated as a $C^*$-algebra by 
$\{i_A(a)i_S(s)\mid a\in A, s\in S\}.$
\end{enumerate} 
\label{def_of_Toeplitz_cp}
\end{definition}
We denote $B$ by $\cT (A, S, \alpha, L)$. If $E$ is a subset of a $C^*$-algebra, 
then $\overline{E}$ is the closed linear span of elements in $E$. 

\begin{definition} A pair $(a, k)$ where $a, k \in \cT (A, S, \alpha, L)$
is a \emph{redundancy} if 
$(a, k)\in i_A(A)\times \overline{i_A(A)i_S(s)
i_S(s)^*i_A(A)}$ for some $s\in S$, and 
\begin{equation}
ai_A(b)i_S(s)=ki_A(b)i_S(s),\text{ for all }b\in A.
\label{def_of_redundancy}
\end{equation}
Let $I$ denote the two-sided 
ideal in $\cT (A, S, \alpha, L)$ generated 
by all possible differences $a-k$, where $(a, k)\in \overline{i_A(A
\alpha_s(A)A)}\times \overline{i_A(A)i_S(s)
i_S(s)^*i_A(A)}$ is a redundancy.

The crossed product is defined to be
\begin{equation}
A\rtimes_{\alpha, L}S:=\cT (A, S, \alpha, L)/I.
\label{def_of_cp}
\end{equation}
\end{definition}

Our immediate goal is to follow \cite[\S3]{F} and construct a product system 
$X$ of Hilbert bimodules over $S$. We then show that the Toeplitz 
algebra $\cT_X$ of $X$ together with appropriate maps $(i_A, i_S)$ satisfies 
the conditions in Definition~\ref{def_of_Toeplitz_cp}, and thus is a Toeplitz 
crossed product for $(A, S, \alpha, L)$. Then we identify the ideal $I$ of 
redundancies in terms of an ideal in $\cT_X$, and accordingly 
describe $A\rtimes_{\alpha, L}S.$

\section{The product system of $(A, S, \alpha, L)$}

\subsection{Preliminaries} 
We recall some definitions and notation from \cite{F}. A Hilbert bimodule 
over a $C^*$-algebra $A$ is a right Hilbert $A$-module $X$ equipped with a 
left action given by a homomorphism $\phi:A\to \mathcal{L}(X)$ 
via the formula 
$a\cdot x:=\phi(a)x$, $a\in A$, $x\in X$. Suppose that $S$ is a countable 
semigroup with 
identity $0$. A family of Hilbert bimodules over $A$ is a \emph{product system} 
over $S$ if there is a semigroup 
homomorphism $p:X\to S$ such that, with $X_s:=p^{-1}(s)$ denoting the fibre 
over $s$, the following condition holds: for every $s, t \in S\setminus\{0\}$, 
the map
$$
(x, y)\in X_s\times Y_t \mapsto xy \in X_{st}
$$
extends to an isomorphism $X_s\otimes_A X_t\cong X_{st}$. Here $X_s\otimes_A 
X_t$ denotes the internal tensor product of Hilbert bimodules, in which elementary  
tensors $x\otimes_A y$ satisfy that $(x\cdot a)\otimes y=x\otimes(a\cdot y)$ 
in $X\otimes Y$, see for example 
\cite[Proposition 3.16]{RW}. Let 
\begin{equation}
\phi_s:A\to \mathcal{L}(X_s), \phi_s(a)(x)=a\cdot x
\label{def_of_left_module_action} 
\end{equation}
denote the homomorphism which defines the left module action. We also 
require that $X_0=A$ with right-action $x\cdot a=xa$ and left action 
$\phi_0(a)(x)=ax$, 
and that the maps $X_0\times X_s\to X_s$, $X_s\times X_0\to X_s$ satisfy 
\begin{equation}
ax=\phi_s(a)x, \text{ and }xa=x\cdot a\text{ for }a\in A, x\in X_s.
\label{properties_of_mult_on_Xe}
\end{equation} 
In order for the map $X_0\times X_s\to X_s$ to induce an isomorphism 
of Hilbert bimodules $X_0\otimes_A X_s\cong X_s$, it is necessary that the bimodule 
$X_s$ be essential, in the sense that 
\begin{equation}
X_s=\clsp\{\phi_s(a)x\mid a\in A, x\in X_s \}.
\label{def_of_essential_module}
\end{equation}

\subsection{The construction of the product system associated with $(A, S, \alpha, 
L)$}.

We assume that a system $(A, S, \alpha, L)$ is given as described in 
Section~\ref{the_system}. 
For each $s\in S$ we endow $A_s:=A\overline{\alpha_s}(1)$ with a right 
$A$-module structure  and a pre-inner product given by
\begin{equation}
x\cdot a:=x\alpha_s(a) \text{ for }a\in A,\text{ and }\langle 
x, y\rangle_s:=L_s(x^*y), 
\forall x,y\in A_s.
\label{right_module_structure_on_Am}
\end{equation}

Notice that when $a$ runs over an approximate unit for $A$ in 
(\ref{def_of_transfer_op}) then $L_s(\overline{\alpha_s}(1)b)=L_s(b)$ 
for all $b\in A$. Thus
\begin{equation}
\langle a\overline{\alpha_s}(1), b\overline{\alpha_s}(1)\rangle_s=L_s(a^*b) 
\label{def_of_pre_inner_product}
\end{equation}
for all $a, b\in A$ and $s\in S$. After modding out elements $x\in A_s$ with 
$\langle x,x\rangle_s=0$ we  let $M_s$ denote the completion of $A_s$ 
to a right Hilbert $A$-module. Since 
$$
\Vert\langle ax, ax\rangle_s \Vert=\Vert L_s(x^*a^*ax)\Vert\leq \Vert 
a^*a\Vert \Vert \langle x, x\rangle_s\Vert
$$
by linearity and positivity of $L_s$, the formula $x\in A_s\mapsto ax\in A_s$ 
extends to an operator $\phi_s(a): 
M_s\to M_s$. A computation using (\ref{def_of_pre_inner_product}) shows 
that $\phi_s(a)$ is adjointable with 
adjoint satisfying $\phi_s(a)^*(y)=a^*y$, and we obtain a homomorphism
$\phi_s:A\to \mathcal{L}(M_s)$. Thus for each $s\in S$,
$X_s:=\{s\}\times M_s$ is a Hilbert bimodule over $A$ with actions
\begin{equation}
(s,x)\cdot a=(s,x\alpha_s(a))\text{ and } \phi_s(a)(s,x)=(s,ax)
\end{equation}
for $a\in A$, $x\in A_s$. We endow $X:=\bigsqcup_{s\in S}X_s$ with 
a semigroup structure as follows. For $x=a\overline{\alpha_s}(1)\in A_s$ and 
$y=b\overline{\alpha_t}(1)\in A_t$ we let
\begin{equation}
(s, x)(t, y)=
(s+t, a\alpha_s(b\overline{\alpha_t}(1)))\in A_{s+t}.
\label{product_on_easy_elements}
\end{equation}
Since $S$ is abelian, $L_{s+t}=L_t\circ L_s$ for all $s, t\in S$, and it follows 
from a straightforward computation that 
$
\Vert a\alpha_s(b\overline{\alpha_{t}}(1))\Vert_{s+t} \leq \Vert 
a\overline{\alpha_s}(1)\Vert_s\Vert b\overline{\alpha_t}(1)\Vert_t. 
$ 
For arbitrary elements $x=\lim_\lambda a_\lambda\overline{\alpha_s}(1)
\in M_s$ and $y=\lim_\lambda b_\lambda\overline{\alpha_t}(1)\in M_t$, this 
estimate shows that $\lim_\lambda a_\lambda
\alpha_s(b_\lambda\overline{\alpha_{t}}(1))$ is a well-defined element of 
$M_{s+t}$, and so we let 
\begin{equation}
(s,x)(t,y):=(s+t, \lim_\lambda a_\lambda
\alpha_s(b_\lambda\overline{\alpha_{t}}(1))).
\label{newproduct_on_easy_elements}
\end{equation}

\begin{prop} The family $X=\bigsqcup_{s\in S}X_s$ of Hilbert bimodules 
over $A$ endowed with the operation 
defined in (\ref{newproduct_on_easy_elements}) is a product system 
over the semigroup $S$. Each fiber $X_s$ is essential 
in the sense of (\ref{def_of_essential_module}).
\label{prodsyst}
\end{prop}

\begin{proof} We let $p(s,x)=s$ for $(s, x)\in X_s$. Using 
(\ref{product_on_easy_elements}) we note that 
(\ref{properties_of_mult_on_Xe}) is immediate when $x\in A_s$, 
and hence is true for all $x\in X_s$. Associativity of multiplication on $X$ 
follows from a routine calculation, which we omit. 
Letting $b$ run over an approximate unit for $A$ in (\ref{product_on_easy_elements}) 
shows that the map $X_s\times X_t\to X_{s+t}$ has dense range, so 
$X_s\otimes_A X_t\to X_{s+t}$ will be an isomorphism  of Hilbert 
bimodules provided that it preserves inner products (see \cite[Lemma 3.2]{F} 
and  \cite[Remark 3.27]{RW}), that is, if
\begin{equation}
\langle\langle u\otimes_A v, u'\otimes_A v'\rangle\rangle
=\langle uv, u'v'\rangle_{s+t},
\label{preserving_inner_products}
\end{equation}
for $u,u'\in X_s$ and $v,v'\in X_t$. If we let 
$u=(s, a\overline{\alpha_s}(1))$, 
$u'=(s, a'\overline{\alpha_s}(1))$, $v=(t, b\overline{\alpha_t}(1))$ and 
$v'=(t, b'\overline{\alpha_t}(1))$ for $a,a',b,b'\in A$,  then the 
left-hand side of 
(\ref{preserving_inner_products}) is
\begin{align}
\langle u', u\rangle_s\cdot v, v'\rangle_t
&=\langle\langle L_s((a')^*a)\cdot v, v'\rangle_t\notag \\
&=\langle (t, L_s((a')^*a))b\overline{\alpha_t}(1), (t, b'
\overline{\alpha_t}(1))\rangle_t\notag \\
&= L_t(b^*L_s(a^*a')b')\notag \\
&=L_t\circ L_s(\alpha_s(b^*)a^*a'\alpha_s(b'))\notag \\
&=L_{s+t}((a\alpha_s(b))^*a'\alpha_s(b'))\text { since }S\text{ is abelian}.\notag
\end{align}
The last expression equals 
$$
\langle(s+t,a\alpha_s(b)\overline{\alpha_{s+t}}(1)), (s+t, a'\alpha_s(b')
\overline{\alpha_{s+t}}(1))\rangle_{s+t},
$$
which is the right-hand side of (\ref{preserving_inner_products}), as 
wanted. Since $\Vert \langle\langle u\otimes_A v, u\otimes_A v\rangle\rangle 
\Vert\leq \Vert u\Vert_s\Vert v\Vert_t$,  
(\ref{preserving_inner_products})  will hold for arbitrary $u\in X_s$ and $v\in 
X_t$. 
If $(u_\lambda)$ denotes an approximate unit for $A$, then for every $b
\in A$, 
$b\overline{\alpha_s}(1)=\lim \phi_s(u_\lambda)b\overline{\alpha_s}(1)$, and 
it follows that $X_s$ is essential.
\end{proof}

\section{The crossed product $A\rtimes_{\alpha, L}S$}

We recall the ingredients necessary to describe the \emph{Toeplitz algebra 
$\cT_X$} from \cite[Proposition 2.8]{F} associated to an arbitrary product 
system $X$ over a semigroup $S$. A map $\psi:X\to B$ 
into a $C^*$-algebra $B$ is called a Toeplitz 
representation of $X$ if $\psi(u)\psi(v)=\psi(uv)$ for all $u,v\in X$, and 
the pair $(\psi_s,\psi_0):=(\psi\vert_{X_s}, \psi\vert_{X_0})$ is a 
Toeplitz representation of the bimodule $X_s$ in the sense of \cite{FR}, i.e.  
$\psi_s$ is linear, $\psi_0:A\to B$ is a homomorphism, and the conditions 
\begin{align}
\psi_s(u\cdot a)&=\psi_s(u)\psi_0(a), \label{Toeplitz_cov_prod_syst1}\\
\psi_s(u)^*\psi_s(v)&=\psi_0(\langle u, v\rangle_s),\text{ and}
\label{Toeplitz_cov_prod_syst2}\\
\psi_s(a\cdot u)&=\psi_0(a)\psi_(u)\label{Toeplitz_cov_prod_syst3}
\end{align}
are satisfied for all $s\in S$, $a\in A$ and $u,v\in X_s$. Note in particular 
that  
$\psi_s$ is bounded. By \cite[Proposition 2.8]{F}, 
there are a $C^*$-algebra $\cT_X$ and a Toeplitz representation $i_X:X 
\to \cT_X$ which enjoy the following universal property: for every Toeplitz 
representation $\psi$ of $X$ there 
is a homomorphism $\psi_*$ on $\cT_X$ for which $\psi_*\circ i_X=\psi$. 
Moreover, $\cT_X$ is generated by $i_X(X)$, and $i_X$ is isometric.

\begin{prop} 
Suppose that $(A, S, \alpha, L)$ is a dynamical system with 
extendible endomorphisms and an action of transfer operators 
as in (\ref{def_of_transfer_op}). Let $X$ be the product system over $S$ 
from Proposition~\ref{prodsyst}, 
let $\cT_X$ be the Toeplitz algebra of $X$ and $i_X:X\to \cT_X$ be the 
universal Toeplitz representation.
\begin{enumerate}
\item[(1).] 
If $(u_\lambda)$ is an approximate unit  for $A$, $(i_X(s, 
\alpha_s(u_\lambda)))_\lambda$ converges strictly in $M(\cT_X)$.
\item[(2).] Let $i_A:A\to \cT_X$ and $i_S:S\to M(\cT_X)$ be given by
$$
i_A(a):=i_X(0,a)\text{ and }i_S(s):=\lim i_X(s, \alpha_s(u_\lambda)) 
$$
for $a\in A$, $s\in S$. Then $(\cT_X, i_A, i_S)$ is a Toeplitz crossed 
product for $(A, S, \alpha, L)$ in the sense of 
Definition~\ref{def_of_Toeplitz_cp}, and $i_A$ is injective.
\end{enumerate}
\label{TX_is_Toeplitz_cp}  
\end{prop}

\begin{proof} Part (1) follows as in \cite[Lemma 3.3]{F}. The rest of the 
proof will follow the argument in the proof of \cite[Proposition 3.4]{F}.  
We verify conditions (i)-(iii) from Definition~\ref{def_of_Toeplitz_cp}. 
Since $i_A(u_\lambda)\to i_X(0, 1_{M(A)})$ for any approximate unit 
$(u_\lambda)$ in $A$, $i_A$ is non-degenerate. Given 
a non-degenerate representation $\sigma$ of $\mathcal{T}_X$, the pair  
$(\sigma \circ i_A, \overline{\sigma}\circ i_S)$ satisfies (TC1)-(TC2), 
because $(i_A, i_S)$ does, as is seen from the calculations
\begin{align}
i_S(s)i_A(a)&=
\lim i_X(s,\alpha_s(u_\lambda)\alpha_s(a))=i_X(s,\overline{\alpha_s}(1)
\alpha_s(a))\notag \\
&=\lim i_X(0,\alpha_s(a))
i_X(s, \alpha_s(u_\lambda))\notag \\
&=i_A(\alpha_s(a))i_S(s), \notag
\end{align}
and 
\begin{align}
i_S(s)^*i_A(a)i_S(s)
&=\lim i_X(s,\alpha_s(u_\lambda))^*i_X(s,a\alpha_s(u_\lambda))  
\notag \\
&=\lim i_X(0, \langle \alpha_s(u_\lambda), a\alpha_s(u_\lambda)\rangle_s)
\text{ by }(\ref{Toeplitz_cov_prod_syst2}) \notag \\
&=\lim i_X(0, L_s(\alpha_s(u_\lambda)^*a\alpha_s(u_\lambda)))\notag \\
&=i_A(L_s(a)).\notag
\end{align}
This proves (i). To prove (ii), suppose that $(\pi, V)$ is a 
Toeplitz covariant representation of $(A,S, \alpha, L)$. Since 
$\Vert\pi(a)V_s \Vert=\Vert L_s(a^*a)\Vert^{1/2}=\Vert a\Vert_s$ by (TC2), 
the equation  
\begin{equation}
\psi(s,a\overline{\alpha_s}(1)):=\pi(a)V_s
\label{def_of_psi_Toeplitz}
\end{equation}
where $s\in S$, $a\overline{\alpha_s}(1)\in A_s$, extends to a well-defined 
map $\psi:X\to B(H_\pi)$, which is multiplicative by property (TC1). 
For $u=(s, 
b\overline{\alpha_s}(1))\in X$ and $a\in A$  we have
\begin{align}
\psi_s(u\cdot a)
&=\pi(b\alpha_s(a))V_s\notag \\
&=\pi(b)V_s\pi(a)\text{ by (TC1) }\notag \\
&=\psi(s, b\overline{\alpha_s}(1))\psi(0,a),\text{ since }V_0=1\notag \\
&=\psi_s(u)\psi_0(a),\notag
\end{align}
giving (\ref{Toeplitz_cov_prod_syst1}), and (\ref{Toeplitz_cov_prod_syst3}) is 
similar. With $u=(s, b\overline{\alpha_s}(1))$  and $v=(s, 
c\overline{\alpha_s}(1))$, 
(\ref{Toeplitz_cov_prod_syst2}) 
follows from 
\begin{align}
\psi_s(u)^*\psi_s(v)
&=V_s^*\pi(b^*c)V_s=\pi(L_s(b^*c))\text{ by (TC2) }\notag \\
&=\psi(0, \langle b\overline{\alpha_s}(1), c\overline{\alpha_s}(1)\rangle_s) 
\text{ since }V_0=1.\notag
\end{align}
Thus $\psi$ is a Toeplitz representation of $X$. By the universal property 
of $\cT_X$ there is a homomorphism $\psi_*:\cT_X\to B(H)$ with $\psi_*\circ 
i_X=\psi$. It is immediate to check that $\pi\times V:=\psi_*$ satisfies the 
first equality in (ii), and the second follows upon evaluating both sides 
on elements in the dense subspace $\{\pi(A)H\}$ of $H$. 
Finally, $i_A(a)i_S(s)=i_X(s, a\overline{\alpha_s}(1))$, and 
these products generate $\cT_X$ because $i_X(X)$ does.
The homomorphism $i_A:A\to \cT_X$ is injective because 
$i_X$ is isometric by \cite[Proposition 2.8]{F}.
\end{proof}

\begin{remark}\label{remark_i_A_i_S_covariant}
Equation (\ref{def_of_psi_Toeplitz}) shows how a Toeplitz representation 
$\psi$ (or $\psi_{\pi, V}$) of $X$ arises from a Toeplitz covariant 
representation $(\pi, V)$ of $(A, S, \alpha, L)$.
Conversely, property (i) in the definition of $\cT(A, S, \alpha, L)$ 
implies that if $\psi:X\to B(H)$ 
is a Toeplitz representation of the associated product system $X$ over 
$S$, then the pair 
$\pi:A\to B(H), V:S\to B(H)$ given by
$$
\pi(a):=(\psi_*\circ i_A)(a)=\psi_0(a)\text{ and } 
V_s=(\overline{\psi_*}\circ i_S)(s),
$$
$a\in A$, $s\in S$, is a Toeplitz covariant representation of 
$(A, S, \alpha, L)$ satisfying  
(\ref{def_of_psi_Toeplitz}).
\end{remark}

Proposition~\ref{TX_is_Toeplitz_cp} shows that $A\rtimes_{\alpha, L}S$ is 
a quotient of $\cT_X$, which we will now identify. Suppose that $\psi:X\to \cT_X$ 
is a Toeplitz representation of $X$, and let $(\psi_s, \psi_0)=(\psi\vert_{X_s}, 
\psi\vert_{X_0})$ be the pair which satisfies 
(\ref{Toeplitz_cov_prod_syst1})-(\ref{Toeplitz_cov_prod_syst3}) for every $s$. Then 
there is a homomorphism $\psi^{(s)}:\mathcal{K}(X_s)\to \cT_X$ such that
$
\psi^{(s)}(\theta_{u,v})=\psi_s(u)\psi_s(v)^*$
and
\begin{equation}
\psi^{(s)}(T)(\psi_s(u))=\psi_s(Tu),
\label{def_of_psi_upper_m}
\end{equation}
for $u,v\in X_s,$
see \cite[Equation (1.1)]{F} and the references therein. 
Let $K_s$ be an ideal in $\phi_s^{-1}(\mathcal{K}(X_s))$ for each $s\in S$,  
and denote by $K$ 
the family of ideals $\{K_s\mid s\in S\}$ in $A$. We define a Toeplitz 
representation  $\psi$ of $X$ to be \emph{coisometric on $K$} if 
$(\psi_s, \psi_0)$ is coisometric on $K_s$ in the sense of 
\cite[Definition 1.1]{FMR} for all $s\in S$, that is if
\begin{equation}
\psi^{(s)}(\phi_s(a))=\psi_0(a)\text{ for all }a\in K_s.
\label{coisometric_rep} 
\end{equation}
If $(\psi_s, \psi_0)$ is coisometric on $\phi_s^{-1}(\mathcal{K}(X_s))$ then 
it is called Cuntz-Pimsner covariant, see \cite{Pi, FR, FMR}, 
and by \cite[Proposition 2.9]{F}, there 
is a $C^*$-algebra, called the Cuntz-Pimsner algebra of $X$, which is universal 
for Toeplitz representations of $X$ such that $(\psi_s, \psi_0)$ is 
Cuntz-Pimsner covariant for all $s\in S$.

\begin{prop} Suppose that $(A, S,\alpha, L)$ is a dynamical 
system with extendible endomorphisms and an action of 
transfer operators as in (\ref{def_of_transfer_op}). Consider the 
product system  $X$ over $S$ from Proposition~\ref{prodsyst}, and the 
universal Toeplitz representation $i_X:X\to \cT_X$. Let  
\begin{equation}
K_s:=\overline{A\alpha_s(A)A}\cap \phi_s^{-1}(\mathcal{K}(X_s))\label{def_of_Ks}
\end{equation}
for all $s\in S$, and denote by $\mathcal{I}_K$  the closed two-sided ideal of 
$\cT_X$ generated by \begin{equation}
\{i_X(0, a)-i^{(s)}(\phi_s(a))\mid s\in S, a\in K_s\}.  
\label{def_of_ideal_in_rel_CP}
\end{equation}

Then $A\rtimes_{\alpha, L}S$ coincides with $\cT_X/\mathcal{I}_K$, and 
the map $j_X:X\to A\rtimes_{\alpha, L}S$ obtained by composing the quotient 
map $\cT_X\to \cT_X/\mathcal{I}_K$ with $i_X$  is 
universal for Toeplitz representations of $X$ which are coisometric on $K=
\{K_s\}_{s\in S}$. 
In particular, if each $\alpha_s$ is surjective or if 
$A$ is unital and $\alpha_s(1)=1$ for all $s\in S$, then $A\rtimes_{\alpha, L}S$ 
is the Cuntz-Pimsner algebra of $X$. 
\label{cp_is_rel_CK_algebra}
\end{prop}

Note that this result motivates that we can 
regard $A\rtimes_{\alpha, L}S$ as the \emph{relative 
Cuntz-Pimsner algebra} of the product system $X$ over $S$ associated to 
$(A, S, \alpha, L)$  and the family of ideals $\{K_s\}_{s\in S}$, where 
$K_s$ is given by equation (\ref{def_of_Ks}).

In general, for an arbitrary product system $X$ over $S$ of Hilbert bimodules over a 
$C^*$-algebra $A$, and for a given family $\{K_s\}_{s\in S}$ of ideals 
$K_s$ in $\phi_s^{-1}(\mathcal{K}(X_s))$, one can mimic the proof of 
\cite[Proposition 2.9]{F} to prove the existence of a $C^*$-algebra 
$\mathcal{O}(X, \{K_s\}_{s\in S})$ and a map $j_X:X\to 
\mathcal{O}(X, \{K_s\}_{s\in S})$  
which are  universal for Toeplitz 
representations of $X$ that are coisometric on $\{K_s\}_{s\in S}$. It seems 
natural to call $\mathcal{O}(X, \{K_s\}_{s\in S})$ the 
\emph{relative Cuntz-Pimsner algebra of $X$ and $\{K_s\}_{s\in S}$}.

\begin{proof}
It suffices by  Proposition~\ref{TX_is_Toeplitz_cp} to show that 
$I=\mathcal{I}_K$ (see also \cite[Lemma 3.5]{BR}). 
 Let $J_s:=\clsp\{i_A(A)i_S(s)\}$ for every $s\in S$. By 
(\ref{def_of_psi_Toeplitz}), 
\begin{align}
i^{(s)}(\theta_{(s, c\overline{\alpha_s}(1)), (s, 
d\overline{\alpha_s}(1))})
&=i_s(c\overline{\alpha_s}(1))
i_s(d\overline{\alpha_s}(1)))^* \notag \\
&=i_A(c)i_S(s)i_S(s)^*i_A(d)^*,\notag
\end{align}
from which it follows by linearity and continuity that 
\begin{equation}
i^{(s)}(\mathcal{K}(X_s))= J_sJ_s^*.
\label{image_of_compacts}
\end{equation}
Let $a\in K_s$. Then 
$(i_X(0, a),i^{(s)}(\phi_s(a)))\in \overline{i_A(A\alpha_s(A)A)}\times 
J_sJ_s^*$ by (\ref{image_of_compacts}), so to show that the pair 
is a redundancy we check for $b\in A$ and $\phi_s(a)=\theta_{u,v}$ with  
$u,v\in X_s$ that
\begin{align}
i^{(s)}(\phi_s(a))i_A(b)i_S(s)
&=i_s(u)i_s(v)^*i_X(s, b\overline{\alpha_s}(1))\text{ by }
(\ref{def_of_psi_Toeplitz})\notag \\
&=i_s(u)i_0(\langle v,b\overline{\alpha_s}(1)\rangle_s)\text{ by }
(\ref{Toeplitz_cov_prod_syst2})\notag\\
&=i_s(u\cdot\langle v,b\overline{\alpha_s}(1)\rangle_s)\text{ by } 
(\ref{Toeplitz_cov_prod_syst1})\notag\\
&=i_s(\theta_{u,v}(b\overline{\alpha_s}(1)))\notag \\
&=i_X(s,ab\overline{\alpha_s}(1))\notag \\
&=i_A(a)i_A(b)i_S(s).\notag 
\end{align}
By linearity and continuity we can extend this identity to all $\phi_s(a)$ 
in $\mathcal{K}(X_s)$, and so it follows that $\mathcal{I}_K\subset I.$

Conversely, suppose that $(i_A(a), k)$ 
is a redundancy in $\overline{i_A(A\alpha_s(A)A)}\times J_sJ_s^*$. Thus 
by (\ref{image_of_compacts}),  
$k=i^{(s)}(T)$ with $T\in \mathcal{K}(X_s)$. Let $x=(s, 
b\overline{\alpha_s}(1))$ in $X_s$ for $b\in A$. Then
\begin{align}
i_s(\phi_s(a)(x))
&=i_X(s, ab\overline{\alpha_s}(1))
=i_A(a)i_A(b)i_S(s)=\notag \\
&=ki_A(b)i_S(s)
=i^{(s)}(T)i_s(x)\text{ by }(\ref{def_of_psi_Toeplitz})\notag \\
&=i_s(Tx)\text{ by }(\ref{def_of_psi_upper_m}),\notag
\end{align}
from which it follows that $\phi_s(a)=T\in \mathcal{K}(X_s)$ because 
$i_s$ is isometric. Hence $a\in K_s$ and  
$I\subset I_K$, and therefore $A\rtimes_{\alpha, L}S=\cT_X/\mathcal{I}_K$. 
To prove that $A\rtimes_{\alpha, L}S$ is universal for Toeplitz representations 
on $X$ which are coisometric on $K$ it suffices by universality of $(\cT_X, i_X)$ 
to show that if $\psi$ is coisometric, then the corresponding homomorphism 
$\psi_*$ on $\cT_X$ vanishes on $\mathcal{I}_K$. This is easily verified on 
elements $i_X(0, a)-i^{(s)}(\phi_s(a))$ when $\phi_s(a)=\theta_{u, v}\in 
\mathcal{K}(X_s)$, whence an $\varepsilon/2$ argument and 
continuity of $\psi_*, i^{(s)}$ and $\psi^{(s)}$ take care of arbitrary 
elements.

When $\alpha_s$ is surjective or $\alpha_s(1)=1$ for all $s$, $K_s=
\phi_s^{-1}(\mathcal{K}(X_s))$, and the proof is complete.
\end{proof}

We know from  \cite[Theorem 4.7]{E} that Exel's  crossed product 
$A\rtimes_{\alpha, L}\N$ of a unital $C^*$-algebra $A$ by a 
single endomorphism $\alpha$ and a transfer operator $L$ is isomorphic 
to $A\rtimes_\alpha \N$ from \cite{ALNR, M}. We recall that 
given an action $\alpha:S\to \End(A)$ of a unital semigroup $S$, the 
\emph{semigroup crossed product} $A\rtimes_\alpha S$ is the universal 
$C^*$-algebra for \emph{covariant} representations  $(\pi, V):(A, S)\to B(H)$, 
where $\pi$ is non-degenerate, $V$ is an isometric representation of $S$ 
and 
$V_s\pi(a)V_s^*=\pi(\alpha_s(a))$, for all $a\in A$, $s\in S$,
see for example \cite{M, LR1} for unital $A$, and   
\cite{Lar} for non-unital $A$. The next result is the generalisation 
of \cite[Theorem 4.7]{E}  to abelian semigroups.

\begin{prop}\label{inj_endom_hereditary_range} 
Suppose that $\alpha:S \to \End(A)$ is an action by extendible injective 
endomorphisms such that $\alpha_s$ has hereditary 
range for all $s\in S$. Let  
\begin{equation}
L_s(a):=\alpha_s^{-1}(\overline{\alpha_s}(1)a\overline{\alpha_s}(1)),\text{ for }
a\in A, s\in S. 
\label{L_from_cond_expectation}
\end{equation}
Then $L$ is an action by transfer operators as in (\ref{def_of_transfer_op}), and 
$A\rtimes_{\alpha, L}S\cong A\rtimes_\alpha S$.  
\end{prop}

\begin{proof}
It follows from \cite[Proposition 4.1]{E} 
that $\alpha_s(A)=\overline{\alpha_s}(1)A\overline{\alpha_s}(1)$ and 
$L_s$ is a transfer operator for $\alpha_s$, for every $s\in S$. Since  
$\alpha_{s+t}=\alpha_t\circ \alpha_s$, a calculation shows that  
$L_{s+t}=L_s\circ L_t$, $\forall s, t\in S$. Let $X$ be the associated 
product system over $S$ of Hilbert bimodules over $A$. We claim that 
\begin{equation}
\phi_s(\overline{A\alpha_s(A)A})\subset \mathcal{K}(X_s) \text{ for all }s\in S.
\label{compacts_is_ideal_alpha_m}
\end{equation}
Indeed, let $a=\sum_{i=1}^n a_i\alpha_s(b_i)c_i^*\in \operatorname{span}\{ 
A\alpha_s(A)A\}$ and $x=(s, b\overline{\alpha_s}(1))\in X_s$. Since 
\begin{align}
\phi_s(a)(x)
&=(s, \sum_{i=1}^n a_i\alpha_s(b_i)c_i^*b\overline{\alpha_s}(1))\notag \\
&=(s, \sum_{i=1}^n a_i\alpha_s(b_i)\alpha_s\circ L_s (c_i^*b))\notag \\
&=\sum_{i=1}^n (s, a_i\alpha_s(b_i))\cdot L_s(c_i^*b)\notag \\
&=\sum_{i=1}^n (s, a_i\alpha_s(b_i))\cdot \langle (s, c_i\overline{\alpha_s}(1)), 
(s, b\overline{\alpha_s}(1))\rangle_s\notag \\
&=\sum_{i=1}^n \theta_{(s, a_i\alpha_s(b_i)), (s, c_i\overline{\alpha_s}(1))}(x),
\label{phi_m(a)_as_compact}
\end{align}
it follows by continuity  that $\phi_s(\overline{A\alpha_s(A)A})\subset 
\mathcal{K}(X_s)$. 
Thus Proposition~\ref{cp_is_rel_CK_algebra} implies that $A\rtimes_{\alpha, L}S$ 
is the universal $C^*$-algebra for coisometric Toeplitz representations of $X$ 
on 
the family $\{K_s:=\overline{A\alpha_s(A)A}\}_s$ of ideals of $A$. The proof of the 
proposition will be concluded by a direct application of the next lemma.
\end{proof}

\begin{lemma}\label{coisometric_as_covariant}
With the hypotheses of  Proposition~\ref{inj_endom_hereditary_range}, let 
$X$ be the product system over $S$ of Hilbert bimodules over $A$ given by 
Proposition~\ref{prodsyst}. Suppose that  $\psi:X\to B(H)$ is a Toeplitz 
representation of $X$, and let $\pi(a)=\psi_0(a), V_s:=\overline{\psi_*}
\circ i_S(s)$, for all $a\in A, s\in S$. We have:
\begin{enumerate}
\item For fixed $s\in S$, $(\psi_s,\psi_0)$ is coisometric on 
$K_s=\overline{A\alpha_s(A)A}$ 
if and only if 
\begin{equation}
\pi(\alpha_s(c))=V_s\pi(c)V_s^*,\forall c\in A.
\label{usual_covariance}
\end{equation} 
\item  $\psi$ is coisometric on $\{K_s\}_{s\in S}$ if 
and only if $(\pi, V)$ is a covariant representation of the triple $(A, S, 
\alpha)$ in the sense of \cite{M, LR1}. 
\end{enumerate}
\end{lemma}

\begin{proof}
By Remark~\ref{remark_i_A_i_S_covariant}, $(\pi, V)$ is a 
Toeplitz covariant representation of $(A, S, \alpha, L)$. Suppose that 
$(\psi_s, \psi_0)$ is coisometric on $K_s$, and 
let $a=b\alpha_s(c)d^* \in K_s$. Then 
\begin{align}
\pi(a)
&=\psi^{(s)}(\phi_s(a))=\psi^{(s)}(\theta_{(s, b\alpha_s(c)), 
(s, d\overline{\alpha_s}(1))})\text{ by }(\ref{phi_m(a)_as_compact}) \notag\\
&=\psi(s, b\alpha_s(c))\psi(s, 
d\overline{\alpha_s}(1))^*\notag 
\\
&=\pi(b\alpha_s(c))V_sV_s^*\pi(d^*)\text{ by }(\ref{def_of_psi_Toeplitz})
\label{pi_coisom_cov}
\end{align}
from which by non-degeneracy of $\pi$ we obtain $
\pi(b\alpha_s(c))=\pi(b\alpha_s(c))V_sV_s^*$. Now letting $b$ run over an 
approximate unit for $A$ and invoking (TC1) gives (\ref{usual_covariance}).

Conversely, (\ref{usual_covariance}) and (TC1) imply 
that $\pi(\alpha_s(c))=\pi(\alpha_s(c))V_sV_s^*$, and we can reverse the 
computations in (\ref{pi_coisom_cov}) and use continuity of $\pi$ 
to conclude that $(\psi_s, \psi_0)$ is coisometric on $K_s$, proving part (i).

Letting $a$ run over an approximate unit for $A$ in 
(TC2) gives that $V_s^*V_s=1$, so $V_s$ is an isometry. 
Part (ii) then follows from (i).
\end{proof}

The next result generalises \cite[Proposition 2.11]{F} 
to relative Cuntz-Pimsner algebras of product systems over the additive semigroup 
$\N$ in the case of the product systems associated to the quadruples
$(A, S, \alpha, L)$ described in Proposition~\ref{inj_endom_hereditary_range}. 
The generalisation of \cite[Proposition 2.11]{F} to 
relative Cuntz-Pimsner algebras of arbitrary product systems over $\N$ 
could go two ways, when either the left action of each fibre is isometric, or 
the left action restricted to a family of ideals $K_n$ of $A$, $n\in \N$, 
is by compacts in each fibre. In the first direction, N. Brownlowe has obtained 
a generalisation in \cite{Br2}.

\begin{prop}\label{rel_CP_as_first_CP} 
Suppose that $\alpha:\N\to \End(A)$ is an action by extendible 
injective endomorphisms 
such that $\alpha_n$ has hereditary range for all $n\in \N$. Let $L$ 
be the action by transfer operators given in (\ref{L_from_cond_expectation}), 
let $X_\N$ be the product 
system over $\N$ of Hilbert bimodules over $A$ associated to $(A, \N, \alpha, 
L)$ by Proposition~\ref{prodsyst}, and let 
$K_n=\overline{A\alpha_n(A)A}$, $n\in \N$. 
Then 
the relative Cuntz-Pimsner algebra of $X_\N$ corresponding to the family of ideals
$(K_n)_n$ is isomorphic to the relative Cuntz-Pimsner algebra $\mathcal{O}(X_1, 
K_1)$, where $X_1$ is the bimodule at the fibre $n=1$.
\end{prop}

\begin{proof}
Let $j_X:X_\N\to A\rtimes_{\alpha, L}\N$ be the universal Toeplitz representation 
of $X_\N$ which is coisometric on $\{K_n\}_n$, and let $(j_{X_1}, j_A):(X_1, A)
\to \mathcal{O}(K_1, X_1)$ be the universal representation of the relative 
Cuntz-Pimsner algebra $\mathcal{O}(K_1, X_1)$, see \cite[Proposition 1.3]{FMR}. 
As in the proof of the similar statement in \cite[Proposition 2.11]{F}, the 
universal property of $j_X$ yields a homomorphism $j_0\times j_1:\mathcal{O}(K_1, 
X_1)\to A\rtimes_{\alpha, L}\N$, and the universal property of $(j_{X_1}, j_A)$ 
combined with \cite[Proposition 1.8(1)]{FR} produces  
a Toeplitz representation $\psi:X_\N\to \mathcal{O}(K_1, X_1)$ such 
that $(\psi_1, \psi_0)=(j_{X_1}, j_A)$. It 
remains to prove that $\psi$ is coisometric on $\{K_n\}_n$, because then $\psi_*$ 
will be the inverse of $j_0\times j_1$, as required. 
 
However, unlike the proof of \cite[Proposition 2.11]{F}, we will 
not use the characterisation of coisometric representations given in 
\cite[Lemma 1.9]{FMR}, since that doesn't seem helpful here. Instead we use 
the relationship between coisometric 
representations of $X_\N$ and representations of $(A, \N, 
\alpha, L)$ as encoded in Lemma~\ref{coisometric_as_covariant}. 

By construction, $(\psi_1, \psi_0)$ is coisometric on $K_1$.  Denote 
$
(\pi, V)=(\psi_0, \overline{\psi_*}\circ i_\N).
$  
Lemma~\ref{coisometric_as_covariant} (i) implies that 
$\pi(\alpha(c))=V_1\pi(c)V_1^*,\forall c\in A. $
Iterating this equation gives $\pi(\alpha^n(c))=V_1^n\pi(c)(V_1^n)^*$ for all 
$n\geq 1$. Since by definition 
$V_n=\lim \psi(n, \alpha^n(u_\lambda))$, where $(u_\lambda)$ is a fixed but arbitrary 
approximate unit for $A$, a computation using multiplicativity of $\psi$ on $X$ 
shows that $V_n=V_1^n$ for all $n\geq 1$. Thus  
$$
\pi(\alpha^n(c))=V_n\pi(c)(V_n)^*,\forall c\in A, \forall n\geq 1,
$$
and so $(\pi, V)$ is a covariant representation of $(A, \N, \alpha)$ as in \cite{M, 
LR1}. Applying Lemma~\ref{coisometric_as_covariant} (ii) finishes the proof, because 
$A\rtimes_{\alpha, L}\N$ is the relative Cuntz-Pimsner algebra of $X_\N$ 
with respect to $\{K_n\}_n.$
\end{proof}

\begin{remark} This result confirms that our definition of $A\rtimes_{\alpha, L}S$ 
is a good extension of Exel's construction from \cite{E}, as we shall now explain. 

Suppose that $(A, S, \alpha, L)$ is a system as in 
Proposition~\ref{inj_endom_hereditary_range}, so that $A\rtimes_{\alpha, 
L}S\cong A\rtimes_\alpha S$. We also assume that $A$ is unital, and we specialise 
to  the case $S=\N$. Thus we infer from \cite[Theorem 4.7]{E} that Exel's  
$A\rtimes_{\alpha, L}\N$ is isomorphic to our  
$A\rtimes_{\alpha, L}\N$. 

The point we make is that Proposition~\ref{rel_CP_as_first_CP} recovers 
Brownlowe and Raeburn's identification in \cite[Proposition 3.10]{BR} of 
Exel's $A\rtimes_{\alpha, L}\N$ with the  relative Cuntz-Pimsner algebra 
$\mathcal{O}(K_1, M_L)$ of a certain Hilbert bimodule $M_L$ over $A$. To see this, 
we need to justify that $X_1$ and $M_L$ represent the same bimodule over $A$. 
Indeed, recall that $M_L$ is the Hilbert module completion of $A/N_L$ in the 
inner product $\langle a, b\rangle_L=L(a^*b)$, where $N_L:=\{a\mid \langle a, 
a\rangle_L=0\}$. In our setting, we have
$$
\{a\alpha(1)\in A_s\mid \langle a\alpha(1), a\alpha(1)\rangle_1=0\}
=\{a\alpha(1)\mid L(\alpha(1)a^*a\alpha(1))=0\}=0,
$$
so $X_1$ is the completion of $A\alpha(1)$ in the inner 
product given by the same formula, and the map $a+N_L\mapsto a\alpha(1)$ 
from $A/N_L$ to $A\alpha(1)$ extends to a module isomorphism of $M_L$ onto $X_1$.
\end{remark}

\begin{remark} With the hypotheses of Proposition~\ref{inj_endom_hereditary_range}, 
suppose that $A$ is commutative and $S=\N$. 
Let $X$ be the product system over $\N$ associated to $(A, \N, \alpha, L)$.
Then the left action in each fibre is by compacts. 
Indeed, simply verify that $\phi(bc^*)=\theta_{b\overline{\alpha}(1), 
c\overline{\alpha}(1)}$ for all $b, c\in A$. Hence \cite[Proposition 2.11]{F}  
implies that the Cuntz-Pimsner algebra of the product system $X$ is canonically   
isomorphic to the Cuntz-Pimsner algebra of the bimodule at fibre $n=1$.   
\end{remark}

\section{Applications and examples}

We now return to our motivating example. We denote by $\N^*$ the multiplicative 
semigroup of non-zero positive integers.

\begin{prop} Let $\Gamma$ be a compact abelian group, and denote by $\varphi$ 
the action of $\N^*$ by endomorphisms 
of $\Gamma$ given by $\varphi_n(s)=s^n$, $\forall n\in \N^*$. Let 
$\beta:\N^*\to \End(C(\Gamma))$ be the action by endomorphisms 
defined by $\beta_n(f)(s):=f(\varphi_n(s))$. Suppose that the following 
three conditions are satisfied by $\varphi$.
\begin{enumerate} 
\item[(I)] $[ \Gamma:\varphi_n(\Gamma) ]<\infty$ for all $n\in \N^*$;
\item[(II)] $\ker \varphi_n$ is finite for all $n\in \N^*$, 
and 
\item[(III)] $\vert \ker\varphi_{mn}\vert= \vert\ker \varphi_{n}
\vert \cdot\vert \ker\varphi_{m}\vert$  
for all $m, n\in \N^*$.  
\end{enumerate}
Then the formula
\begin{equation}
L_n(f)(t):=\begin{cases}\frac 1{\vert \{s\in \Gamma\mid \varphi_n(s)=t
    \}\vert}\sum_{\{s\in \Gamma\mid \varphi_n(s)=t\}}f(s)&\text{ if }t\in 
\varphi_n(\Gamma)\\
0&\text{ otherwise}\end{cases}
\label{def_of_L}
\end{equation}
defines an action by transfer operators $L_n$ for $\beta_n$, $\forall n\in N^*$.
\label{beta_L_when_injective}
\end{prop}

\begin{proof}
Note that $\varphi_n(\Gamma)$ is a closed subgroup of $\Gamma$ which is also 
open because it has finite index. Thus $L_n$ is well-defined for every $n\in 
\N^*$. A direct computation shows that $L_n(f\beta_n(g))=L_n(f)g$ for all 
$f, g\in C(\Gamma)$, so $L_n$ is a transfer operator for $\beta_n$. 
 To see that $L$ defines an action, we compute that 
\begin{align}
L_m(L_n(f))(t)&=\begin{cases}\frac 1{\vert \{s\mid s^m=t\}\vert}\frac 1{\vert 
\{r\mid r^n=s\}}\vert \sum_{\underset{\{r\mid r^n=s\}}{\{s\mid s^m=t\}}}
f(r)&\text{ if }t\in 
\varphi_m(\Gamma), s\in \varphi_n(\Gamma)\\
0&\text{ otherwise}\end{cases}\notag \\
&=\begin{cases}\frac 1{\vert \{r\mid r^{mn}=t\}\vert}\sum_{\{r\mid r^{mn}=t\}}
f(r)&\text{ if }t\in 
\varphi_{mn}(\Gamma)\\
0&\text{ otherwise}\end{cases}\notag \\
&=L_{mn}(f)(t),\label{L_mult}
\end{align}
because 
$
\{r\mid r^{mn}=t\}=\bigcup_{\{s\mid s^m=t\}}\{r\mid r^n=s\}
$
and $\vert\{r\mid r^{mn}=t\}\vert=\vert  \{r\mid r^n=s\}\vert 
\cdot\vert \{s\mid s^m=t\} \vert$ by (III).
\end{proof}

It is important to stress that the crossed product $C(\Gamma)
\rtimes_{\beta, L}\N^*$ 
is a new example, as the theory we have so far does not obviously identify 
it with an algebra of the form $A\rtimes_\alpha \N^*$ for some action $\alpha$. 
We will be able to say more if we concentrate attention to two special cases. 
One arises if we  relax condition (I) to say that 
$\varphi_n$ is surjective for all $n\in \N^*$, the other if we ask in (II) that 
$\varphi_n$ is injective for all $n\in \N^*$. 

\subsection{The case when all $\varphi_n$ are surjective}\label{all_phi_surj}

We first observe that if $\varphi_n$ is a surjective homomorphism  
for all $n\in \N^*$, then the map $z\mapsto \varphi_n(z)$ is a homomorphism 
of $\ker \varphi_{mn}$ onto $\ker\varphi_m$ with kernel equal to 
$\ker\varphi_n$. Thus, if in addition condition (II) holds, then (III) is also 
satisfied, and Proposition~\ref{beta_L_when_injective} takes the following 
form.

\begin{prop}\label{all_phin_surj}
Let $\Gamma$ be a compact abelian group, and suppose that $\varphi_n(s)=s^n$ for 
$n\in \N^*$ defines an action of $\N^*$ by surjective endomorphisms of $\Gamma$ 
with finite kernel. Then $\beta_n(f)(s)=f(\varphi_n(s))$ implements 
an action by (unital) injective endomorphisms of $C(\Gamma)$, and 
$$
L_n(f)(t)=\frac 1{\vert \{s\mid \varphi_n(s)=t \}\vert}\sum_{\{s\in \Gamma\mid 
\varphi_n(s)=t\}}f(s)
$$
for $n\in \N^*$, $f\in C(\Gamma)$, is an action by transfer operators for $\beta$. 
\end{prop}

As an example we can take for $\Gamma$ the compact divisible group $\T$. 
It would 
be interesting to identify $C(\T)\rtimes_{\beta, L}\N^*$ in terms of familiar 
$C^*$-algebras. Since $\beta_n(1)=1$, $\forall n\in \N^*$, the crossed product 
by $\beta$ in the 
sense of \cite{M, LR1} is not meaningful here, because the covariance relation 
$V_n \pi(a) V_n^*=\pi(\beta_n(b))$ would imply that $V_n$ is a unitary.

\subsection{The case when all $\varphi_n$ are injective.}

We now relax condition (II) to say that $\varphi_n$ is injective for all $n\in 
\N^*$. Note in particular that (III) is trivially satisfied. 
Proposition~\ref{beta_L_when_injective} now takes the form of the following result, 
whose straightforward proof we omit.

\begin{prop}\label{general_form_varphi_inj}
 Let $\Gamma$ be a compact abelian group, let $\varphi$ be the action 
of $\N^*$ by endomorphisms of $\Gamma$ given by $\varphi_n(s)=s^n$ for $n\in \N^*$, 
and let $\beta:\N^*\to \End (C(\Gamma))$ be given by $\beta_n(f)(t)=f(\varphi_n(t))$ 
for $n\in \N^*$ and $f\in C(\Gamma)$. 

If $\varphi_n$ is injective for all $n\in \N^*$ and if condition (I) holds, then 
$\beta$ is an action by surjective endomorphisms, (\ref{def_of_L}) becomes 
$$
L_n(f)(t)=\begin{cases}f(\varphi_n^{-1}(t))&\text{ if }t\in 
\varphi_n(\Gamma)\\
0&\text{ otherwise},\end{cases}
$$
and $L$ is an action of $\N^*$ by endomorphisms of $C(\Gamma)$.  
\end{prop}

The next proposition contains a general construction of 
systems $(C(\Gamma), \N^*, \beta, L)$ of the form described in 
Proposition~\ref{general_form_varphi_inj}.

\begin{prop}\label{general_stuff}
Let $G$ be a discrete abelian group, and let $\phi$ be the action of $\N^*$ 
by endomorphisms of $G$ given by  $\phi_n(x)=nx$, for $x\in G$ and $n\in \N^*$. 

Suppose that $\phi_n$ is surjective and $\ker \phi_n$ is finite, $\forall n\in 
\N^*$. Then  the following hold. 
\begin{enumerate}
\item The homomorphism
$\widehat{\phi_n}:\widehat{G}\to\widehat{G}$ given by
$\chi\mapsto \chi\circ \phi_n$ is injective and satisfies that
$[ \widehat{G} :\widehat{\phi_n}(\widehat{G})]
=\vert \operatorname{ker}\phi_n\vert <\infty$, for all $n\in \N^*$.
\item There is an action $\alpha:\N^*\to \End(C^*(G))$ given by 
\begin{equation}
\alpha_n(\delta_y)=\frac 1{\vert
  \operatorname{ker}\phi_n\vert}\sum_{\{x\in G\mid \phi_n(x)=y\}}\delta_x, \forall 
n\in \N^*. 
\label{formula_action_on_group_alg}
\end{equation}
\item The Fourier transform action of $\alpha$ on $C(\widehat{G})$ is given by
\begin{equation}
\widehat{\alpha_n}(f)(\chi)=\begin{cases}
  f({\widehat{\phi_n}}^{-1}(\chi)),&\text{ if
  }\chi\in\widehat{\phi_n}(\widehat{G})\\
0,&\text{ if }\chi\notin\widehat{\phi_n}(\widehat{G}).
\end{cases}
\label{formula_action_on_cont_funct}
\end{equation}
\end{enumerate}
\end{prop}

Before proving this we make a note of an immediate consequence and we discuss 
the connection with \cite{LarR}.

\begin{cor}\label{system_dual_system}
With the notation and assumptions of Proposition~\ref{general_stuff}, 
let $\widehat{\beta_n}(f)(\chi)=f(\widehat{\phi_n}(\chi))$ for $f\in 
C(\widehat{G})$, $\chi \in \widehat{G}$, $n\in \N^*$. Then $(C(\widehat{G}), 
\N^*, \widehat{\beta}, \widehat{\alpha})$ is an example of a system as in 
Proposition~\ref{general_form_varphi_inj}. 

Dually, $\alpha$ from Proposition~\ref{general_stuff} (ii) is an action of $\N^*$ 
by transfer operators for $\beta:\N^*\to \End(C^*(G))$ defined by $\beta_n(\delta_x)
=\delta_{nx}$. 
\end{cor}

\begin{remark}
It seems interesting to recall the construction of endomorphic actions of  
$\N^*$ (indeed, of the additive semigroup $\N^k$ for $k\in \{1, 2, \dots\}\cup 
\{\infty\}$, where $\N^*$ is identified as $\N^\infty$) from 
\cite{LarR}, and to discuss its possible connection with 
Proposition~\ref{general_stuff}. In the setup of \cite[\S 1]{LarR}, one 
considers an action $\eta$ of a discrete abelian group $G$ by injective 
endomorphisms whose range $\eta_m(G)$ is a subgroup of $G$ of finite index 
for all $m\in \N^*$. Then \cite[Proposition 1.3]{LarR} establishes the existence 
of an action of $\N^*$ by endomorphisms of $C^*(G_\infty/G)$, where $G_\infty$ 
is a certain direct limit built from $G$ and $\eta$. In the case when 
$\eta_m(g)=mg$ for $m\in \N^*$ and $g\in G$, Proposition~\ref{general_stuff} 
explains why the construction of \cite{LarR} was possible: indeed, under the  
given assumptions $\eta$ induces an action of $\N^*$ by surjective 
endomorphisms of $G_\infty/G$ with finite kernel. Therefore equation 
(1.5) of \cite{LarR} follows directly from  Proposition~\ref{general_stuff} (ii). 
\end{remark}

\begin{proof}[Proof of Proposition~\ref{general_stuff}.] 
Fix $n\in \N^*$ and suppose $\widehat{\phi_n}(\chi)=1$. Then
$\chi(\phi_n(x))=1$ for all $x\in G$, and hence $\chi(y)=1$ for all $y\in
G$ because $\phi_n$ is onto. Thus $\widehat{\phi_n}$ is injective. 

Since $\phi_n$ is an isomorphism of $G/\ker \phi_n$ onto $G$, 
$\widehat{\phi_n}$ is an isomorphism of $\widehat{G}$ onto $(G/\ker 
\phi_n)^\wedge$. Therefore $\widehat{\phi_n}(\widehat{G})=(G/\ker 
\phi_n)^\wedge\cong (\ker \phi_n)^\perp$. Because $\chi
\mapsto \chi\vert_{\operatorname{ker}\phi_n}$ is a surjection of
$\widehat{G}$ onto $\widehat{\operatorname{ker}\phi_n}$ with 
kernel $(\operatorname{ker}\phi_n)^\perp$, it follows that 
$\widehat{G}/\widehat{\phi_n}(\widehat{G})=
\widehat{G}/(\operatorname{ker}\phi_n)^\perp\cong
\widehat{\operatorname{ker}\phi_n }
$, proving claim (i). 

To prove (ii), we first show that the Fourier transform of the right hand side of
(\ref{formula_action_on_group_alg}) is given by the formula in the right hand 
side of (\ref{formula_action_on_cont_funct}). 
Since (\ref{formula_action_on_cont_funct}) is a well-defined 
endomorphism of $C(\widehat{G})$ because $\widehat{\phi_n}(\widehat{G})$ 
is clopen by (i), we conclude by taking the inverse 
Fourier transform that 
(\ref{formula_action_on_group_alg}) is an endomorphism of $C^*(G)$. This also 
proves (iii).

Take $z\in G$ such that
$\phi_n(z)=y$. Then $\{x\mid \phi_n(x)=y\}=z(\operatorname{ker}\phi_n)$, and for
$\gamma\in \widehat{G}$ we have
\begin{align}
\bigl(\sum_{\phi_n(x)=y}\delta_x \bigr)^\wedge(\gamma)
&=\sum_{\phi_n(x)=\phi_n(z)}\gamma(x)\notag \\
&= \gamma(z)\sum_{\{x\in \ker \phi_n\}}\gamma(x)\notag \\
&=\begin{cases}
\gamma(z)\vert \ker \phi_n\vert,&\text{ if }\gamma\in (\ker \phi_n)^\perp\\
0,&\text{ if }\gamma\notin (\ker \phi_n)^\perp,\notag
\end{cases}
\end{align}
because $\{\gamma(x)\mid x\in \operatorname{ker}\phi_n\}$ is a non-trivial
subgroup of $\T$ when $\gamma \notin  (\operatorname{ker}\phi_n)^\perp$. A routine 
calculation shows that $\gamma(z)=\widehat{\delta_{\phi_n(z)}}
(\widehat{\phi_n}^{-1}(\gamma))$ when $\gamma \in  \widehat{\phi_n}(\widehat{G})=
(\operatorname{ker}\phi_n)^\perp$. So
$$
\bigl(\frac 1{\vert\ker \phi_n\vert}\sum_{\phi_n(x)=\phi_n(z)}
\delta_x \bigr)^\wedge(\gamma)
=\begin{cases}
\widehat{\delta_{\phi_n(z)}}
(\widehat{\phi_n}^{-1}(\gamma)),&\text{ if }\gamma\in \widehat{\phi_n}(\widehat{G})\\
0,&\text{ if }\gamma\notin \widehat{\phi_n}(\widehat{G}),
\end{cases}
$$
as claimed.
To prove that $\alpha_m\alpha_n=\alpha_{mn}$ for $m, n\in \N^*$ we 
need to know that $\vert \ker \phi_{mn}\vert=\vert \ker \phi_m\vert \cdot \vert 
\ker \phi_n\vert$, and this follows from the observation at the beginning of 
section~\ref{all_phi_surj}. 
\end{proof}

Towards identifying the crossed product $C^*(G)\rtimes_{\beta,\alpha}\N^*$ 
of the system obtained in Corollary~\ref{system_dual_system} we first prove a 
general result, and then notice that we can apply it in our setting because 
$\beta_n\circ\alpha_n=\operatorname{id}$ for all $n\in \N^*$. We stress that 
the setup of the next result combines two actions with different properties: one 
is by unital endomorphisms, for which a crossed product as in \cite{M, LR1} 
implemented by an isometric representation of $S$ is not meaningful, the other 
action is by injective endomorphisms with hereditary range, and this class 
fits very well the crossed product construction of \cite{M, LR1}. The point is 
that by viewing the two actions as transfer operators for each other we can treat 
them in a unified way using our generalisation of Exel's crossed product.

\begin{theo}\label{swapping_Exel_cp}
Let $A$ be a unital $C^*$-algebra. Suppose that $\beta, L:S\to \End A$ are 
actions such that  $\beta_s(1)=1$, $L_s$ is a 
transfer operator for $\beta_s$,    
and $\beta_s\circ L_s=\operatorname{id}$, so that in particular $\beta_s$ 
is trivially a transfer operator for $L_s$, for all $s\in S$.

Then the map $(\pi, V)\mapsto (\pi, V^*)$ is a bijective correspondence 
between Toeplitz covariant representations of $(A, S, L, \beta)$ which satisfy 
\begin{equation}
V_s\pi(a)V_s^*=\pi(L_s(a)), \label{pi_V_cov}
\end{equation}
for all $s\in S, a\in A$, and Toeplitz 
covariant representations  of $(A, S, \beta, L)$ such that 
\begin{equation}
\pi(a)=\pi(a)V_s^*V_s, \label{pi_V*_cov}
\end{equation} 
for all $s\in S, a\in A$, and hence induces an
isomorphism of $A\rtimes_{L, \beta}S$ onto $A\rtimes_{\beta, L}S.$ 
\end{theo}

\begin{proof} 
Denote by $X$ and $Y$ the product systems over $S$ of Hilbert bimodules over $A$ 
associated to $(A, S, \beta, L)$  and $(A, S, L, \beta)$. 

Since $L_s$ is a transfer operator for $\beta_s$, we have 
$L_s\circ \beta_s(b)=L_s(1)bL_s(1)$ for all $s\in S, b\in A$. Thus, since 
$\beta_s$ is surjective, the range of $L_s$ is the corner $L_s(1)AL_s(1)$ of 
$A$ and  $\beta_s(a)=L_s^{-1}(L_s(1)aL_s(1))$, for all $a\in A$, $s\in S$. 
By the proof of Proposition~\ref{inj_endom_hereditary_range},
$A\rtimes_{L, \beta}S$ is the universal $C^*$-algebra for Toeplitz 
representations of $Y$ which are coisometric with respect to 
the family of ideals $\{K_s=\overline{AL_s(A)A}\subset \phi_s^{-1}(\mathcal{K}
(Y_s))\mid s\in S\}$. Since $\beta_s(1)=1$ for all $s$,  
Proposition~\ref{cp_is_rel_CK_algebra} says that  $A\rtimes_{\beta, L}S$ 
is the Cuntz-Pimsner algebra of $X$, and is therefore universal 
for Toeplitz representations of $X$ which are coisometric on the family 
of ideals $\{J_s=\phi_s^{-1}(\mathcal{K}(X_s))\mid s\in S\}$.

Suppose first that $(\pi, V)$ is a Toeplitz covariant representation of 
$(A, S, L, \beta)$ which satisfies (\ref{pi_V_cov}). Since $\beta_s(1)=1$, (TC2) 
implies that $V_s^*V_s=\pi(\beta_s(1))=1$, so $V_s$ is an isometry, hence 
(\ref{pi_V*_cov}) is trivial. Since
\begin{align}
\pi(\beta_s(a))V_s^*
&=V_s^*\pi(L_s\circ \beta_s(a)) \text{ by (TC1) for }(\pi, V)\notag \\
&=V_s^*\pi(L_s(1)a)\notag \\
&=V_s^*\pi(a) \text{ by (\ref{pi_V_cov})},\notag
\end{align}
it follows that $(\pi, V^*)$ is Toeplitz covariant for $(A, S, \beta, L)$ 
because condition (TC2) is identical to (\ref{pi_V_cov}).

Conversely, suppose that $(\pi, V^*)$ is Toeplitz covariant for $(A, S, \beta, 
L)$. Then (TC2) gives (\ref{pi_V_cov}). Using (\ref{pi_V*_cov}) and (TC2) gives 
$
V_s\pi(a)=V_s\pi(a)V_s^*V_s=\pi(L_s(a))V_s$, 
and by (TC1) and (\ref{pi_V*_cov}) we also have
$
V_s^*\pi(a)V_s=\pi(\beta_s(a))V_s^*V_s=\pi(\beta_s(a)),
$
for all $a\in A$ and $s\in S$. Therefore $(\pi, V)$ is a Toeplitz covariant 
representation 
of $(A, S, L,\beta)$ with the required properties.

For every Toeplitz representation $(\pi, V)$ of $(A, S, L, \beta)$ which 
satisfies 
(\ref{pi_V_cov}), let $\psi_{\pi, V}$ denote the Toeplitz representation of 
$Y$ given 
by (\ref{def_of_psi_Toeplitz}) for $(\pi, V)$. By the first part 
of the proof, $\psi_{\pi, V^*}$ is a Toeplitz representation of $X$. 
Lemma~\ref{coisometric_as_covariant} 
implies that $\psi_{\pi, V}$ is coisometric on $\{K_s\}_s$, and we claim that 
$\psi_{\pi, V^*}$ is coisometric on $\{J_s\}_s$. 
Indeed, for  $a\in J_s$ written as $a=bc^*$ for $b, c\in A$, the claim 
follows from the computations 
\begin{align}
\psi_{\pi, V^*}^{(s)}(\phi_s(a))&=\psi_{\pi, V^*}^{(s)}(\theta_{(s, b), (s, 
c)})\notag \\
&=\psi_{\pi, V^*}(s, b)\psi_{\pi, V^*}(s, c)^*=\pi(b)V_s^*V_s\pi(c^*)\notag\\
&=\pi(a) \text{ by }(\ref{pi_V*_cov}).\notag
\end{align}

Thus the bijective correspondence $(\pi, V)\mapsto (\pi, V^*)$ 
induces a bijective correspondence $\psi_{\pi, V}\mapsto \psi_{\pi, V^*}$ from 
the Toeplitz representations of $Y$ which are coisometric on $\{K_s\}_s$ onto 
the Toeplitz representations of $X$ which are coisometric on $\{J_s\}_s$. Hence the 
isomorphism $A\rtimes_{L,\beta} S\cong A\rtimes_{\beta, L}S$ follows from  
universal properties.
\end{proof}

\begin{cor}\label{three_cp_are_iso}
 Under the assumptions of Theorem~\ref{swapping_Exel_cp},
  there are isomorphisms
$$
A\rtimes_{\beta, L}S \cong A\rtimes_{L, \beta}S \cong
A\rtimes_{L}S. 
$$
\end{cor}

\begin{proof} The second isomorphism follows from 
Proposition~\ref{inj_endom_hereditary_range}.
\end{proof}

Notice that the above isomorphisms convert the isometric representation of 
$S$ implementing $A\rtimes_\alpha S$ into a representation of $S$ by coisometries 
which implements $A\rtimes_{\beta, L}S$.

\begin{example}\label{the_BC_algebra} 
The hypotheses of Proposition~\ref{general_stuff} 
are satisfied by the divisible group $\Q/\Z$, since $\vert  \ker \phi_n\vert=n$ 
for all $n\in \N^*$. The crossed product $C((\Q/\Z)^\wedge)
\rtimes_{\widehat{\beta},\widehat{\alpha}}\N^*$, or dually, $C^*(\Q/\Z)
\rtimes_{\beta, 
\alpha}\N^*$, is by Corollary~\ref{three_cp_are_iso} isomorphic to 
$C^*(\Q/\Z)\rtimes_{\alpha}\N^*$. By \cite[Corollary 2.10]{LR2}, this 
last crossed product is isomorphic to the Hecke $C^*$-algebra of Bost and Connes 
from \cite{BC}.
\end{example}

\begin{remark} In Proposition~\ref{beta_L_when_injective} we can restrict 
attention  to the unital subsemigroup $\N^F$ of  $\N^*$ consisting of non-zero  
positive integers which have prime factorisation in a fixed 
subset  $F$ of the set $\mathcal{P}$ of prime numbers. Then it is possible to 
translate the results of this section to the similar actions $\beta$ and $L$ 
of $\N^F$. In particular, the analogue 
of Example~\ref{the_BC_algebra} would show that $C(\Q/\Z)\rtimes_{\beta, 
L}\N^{F}$ is isomorphic to one of the Hecke $C^*$-algebras 
from \cite{BLPR}.
\end{remark}

\end{document}